\documentclass[11pt,reqno]{amsart}
\usepackage{amsmath, amssymb, amsthm}

 \renewcommand{\a}{\alpha}
\renewcommand{\b}{\beta}

\renewcommand{\(}{\left\(}
\renewcommand{\)}{\right\)}
\renewcommand{\[}{\left\[}
\renewcommand{\]}{\right\]}
\renewcommand{\i}{\infty}
\numberwithin{equation}{section}
 \theoremstyle{plain}
\newtheorem{theorem}{Theorem}[section]

   \makeatletter
\def\proof{\@ifnextchar[{\@oproof}{\@nproof}}
\def\@oproof[#1][#2]{\trivlist\item[\hskip\labelsep\textit{#2 Proof of\
#1.}~]\ignorespaces}
\def\@nproof{\trivlist\item[\hskip\labelsep\textit{Proof.}~]\ignorespaces}

\makeatother

\begin{document}
\title[A simple proof of a congruence]{A simple proof of a congruence for a series involving the little $q$-Jacobi polynomials}

\author{Atul Dixit}
\address{Department of Mathematics, Indian Institute of Technology Gandhinagar, Palaj, Gandhinagar, Gujarat 382355, India}\email{adixit@iitgn.ac.in}



\dedicatory{Dedicated to Professor George E. Andrews on the occasion of his 80th birthday}
\begin{abstract}
We give a simple and a more explicit proof of a mod $4$ congruence for a series involving the little $q$-Jacobi polynomials which arose in a recent study of a certain restricted overpartition function. 
\end{abstract}
\maketitle
\section{Introduction}\label{intro}
In \cite{adsy1}, Andrews, Schultz, Yee and the author studied the overpartition function $\overline{p}_{\omega}(n)$, namely, the number of overpartitions of $n$ such that all odd parts are less than twice the smallest part, and in which the smallest part is always overlined. In the same paper, they obtained a representation for the generating function of $\overline{p}_{\omega}(n)$ in terms of a ${}_3\phi_{2}$ basic hypergeometric series and an infinite series involving the little $q$-Jacobi polynomials. The latter are given by \cite[Equation (3.1)]{aa}
\begin{equation}\label{little}
p_n(x;\a,\b:q):={}_{2}\phi_{1}\bigg(\begin{matrix}q^{-n},& \a\b q^{n+1}\\
&\a q\end{matrix}\, ;qx\bigg),
\end{equation}
where the basic hypergeometric series ${}_{r+1}\phi_{r}$ is defined by
\begin{equation*}\label{bhs}
{}_{r+1}\phi_{r}\left(\begin{matrix} a_1, a_2, \ldots, a_{r+1}\\
  b_1,  b_2, \ldots, b_{r} \end{matrix}\,; q,
z \right) :=\sum_{n=0}^{\infty} \frac{(a_1;q)_n (a_2;q)_n \cdots (a_{r+1};q)_n}{(q;q)_n (b_1;q)_n \cdots (b_{r};q)_n} z^n,
\end{equation*}
and where we use the notation
\begin{align*}
&(A;q)_0 =1;\quad (A;q)_n  = (1-A)(1-Aq)\cdots(1-Aq^{n-1}),\quad n \geq 1,\\
&(A;q)_{\i}  = \lim_{n\to\i}(A;q)_n\hspace{3mm} (|q|<1)\nonumber.
\end{align*}
The precise representation for the generating function of $\overline{p}_{\omega}(n)$ obtained in \cite{adsy1} is as follows.
\begin{theorem}\label{fgfthm}
The following identity holds for $|q|<1$:
\begin{align}\label{fgf}
\overline{P}_{\omega}(q):=\sum_{n=1}^{\infty}\overline{p}_{\omega}(n)q^n&=-\frac{1}{2}\frac{(q;q)_{\infty}(q;q^2)_{\infty}}{(-q;q)_{\infty}(-q;q^2)_{\infty}}{}_{3}\phi_{2}\left(\begin{matrix}&-1, &iq^{1/2},&-iq^{1/2}\\
& q^{1/2},& -q^{1/2}
\end{matrix}\, ;q, q \right)\nonumber\\
&\quad+\frac{(-q;q)_{\infty}}{(q;q)_{\infty}}\sum_{n=0}^{\infty}\frac{(q;q^2)_n(-q)^n}{(-q;q^2)_{n}(1+q^{2n})}p_{2n}(-1;q^{-2n-1},-1:q).
\end{align}
\end{theorem}
Later, Bringmann, Jennings-Shaffer and Mahlburg \cite[Theorem 1.1]{bjm} showed that $\overline{P}_{\omega}(q)+\frac{1}{4}-\frac{\eta(4\tau)}{2\eta(2\tau)^2}$, where $q=e^{2\pi i\tau}$ and $\eta(\tau)$ is the Dedekind eta function, can be completed to a function $\hat{P}_{\omega}(\tau)$, which transforms like a weight $1$ modular form. They called the function $\overline{P}_{\omega}(q)+\frac{1}{4}-\frac{\eta(4\tau)}{2\eta(2\tau)^2}$ a \emph{higher depth mock modular form}.

While the series involving the little $q$-Jacobi polynomials in Theorem \ref{fgfthm} itself looks formidable, it was shown in \cite[Theorem 1.3]{adsy1} that modulo $4$, it is a simple $q$-product. The mod $4$ congruence proved in there is given below.
\begin{theorem}\label{congmod4}
The following congruence holds:
\begin{equation}\label{mod4}
\sum_{n=0}^{\infty}\frac{(q;q^2)_n(-q)^n}{(-q;q^2)_{n}(1+q^{2n})}p_{2n}(-1;q^{-2n-1},-1:q)\equiv\frac{1}{2} \frac{(q;q^2)_{\infty}}{(-q;q^2)_{\infty}} \pmod{4}.
\end{equation}
\end{theorem}
The proof of this congruence in \cite{adsy1} is beautiful but somewhat involved. The objective of this short note is to give a very simple proof of it. In fact, we derive it as a trivial corollary of the following result.

\begin{theorem}\label{congmod41}
For $|q|<1$, we have
\begin{align}
&\sum_{n=0}^{\infty}\frac{(q;q^2)_n(-q)^n}{(-q;q^2)_{n}(1+q^{2n})}p_{2n}(-1;q^{-2n-1},-1:q)\nonumber\\
&=\frac{1}{2} \frac{(q;q^2)_{\infty}}{(-q;q^2)_{\infty}}+\frac{4q^2}{(1+q)}\sum_{n=0}^{\infty}\frac{(q^3;q^2)_n(-q)^n}{(-q^3;q^2)_{n}(1+q^{2n+2})}\sum_{j=0}^{n}\frac{(-q;q)_{2j}q^{2j}}{(q^2;q)_{2j}}.
\end{align}
\end{theorem}
The presence of $4$ in front of the series on the right-hand side in the above equation immediately implies that Theorem \ref{congmod4} holds.
\section{Proof of Theorem \ref{congmod41}}
Observe that from \eqref{little},
\begin{align}\label{diff1}
\sum_{n=0}^{\infty}\frac{(q;q^2)_n(-q)^n}{(-q;q^2)_{n}(1+q^{2n})}p_{2n}(-1;q^{-2n-1},-1:q)=
\sum_{n=0}^{\infty}\frac{(q;q^2)_n(-q)^n}{(-q;q^2)_{n}(1+q^{2n})}\sum_{j=0}^{2n}\frac{(-1;q)_j}{(q;q)_j}(-q)^j.
\end{align}
However, let us first consider 
\begin{equation}\label{diff2}
A(q):=\sum_{n=0}^{\infty}\frac{(q;q^2)_n(-q)^n}{(-q;q^2)_{n}(1+q^{2n})}\sum_{j=0}^{2n}\frac{(-1;q)_j}{(q;q)_j}q^j.
\end{equation}
The only difference in the series on the right-hand side of \eqref{diff1} and the series in \eqref{diff2} is the presence of $(-1)^{j}$ inside the finite sum in the former. 

To simplify $A(q)$, we start with a result of Alladi \cite[p.~215, Equation (2.6)]{alladi}:
\begin{equation}\label{alla}
\frac{(abq;q)_n}{(bq;q)_n}=1+b(1-a)\sum_{j=1}^{n}\frac{(abq;q)_{j-1}q^j}{(bq;q)_j}.
\end{equation}
Let $a=-1, b=1$ and replace $n$ by $2n$ so that
\begin{equation}\label{aftall}
\sum_{j=0}^{2n}\frac{(-1;q)_jq^j}{(q;q)_j}=\frac{(-q;q)_{2n}}{(q;q)_{2n}}.
\end{equation}
Substitute \eqref{aftall} in \eqref{diff2} to see that
\begin{align}\label{diff3}
A(q)&=\sum_{n=0}^{\infty}\frac{(q;q^2)_n(-q)^n}{(-q;q^2)_{n}(1+q^{2n})}\frac{(-q;q)_{2n}}{(q;q)_{2n}}\nonumber\\
&=\frac{1}{2}+\sum_{n=1}^{\infty}\frac{(-q^2;q^2)_{n-1}}{(q^2;q^2)_n}(-q)^n\nonumber\\
&=\frac{1}{2}\frac{(q;q^2)_{\infty}}{(-q;q^2)_{\infty}},
\end{align}
where in the last step we used the $q$-binomial theorem $\sum_{n=0}^{\infty}\frac{(a;q)_n}{(q;q)_n}z^n=\frac{(az;q)_{\infty}}{(z;q)_{\infty}}$, valid for $|z|<1$ and $|q|<1$.

From \eqref{diff1} and \eqref{diff2},
{\allowdisplaybreaks\begin{align}\label{diff4}
&\sum_{n=0}^{\infty}\frac{(q;q^2)_n(-q)^n}{(-q;q^2)_{n}(1+q^{2n})}p_{2n}(-1;q^{-2n-1},-1:q)-A(q)\nonumber\\
&=\sum_{n=0}^{\infty}\frac{(q;q^2)_n(-q)^n}{(-q;q^2)_{n}(1+q^{2n})}\sum_{j=0}^{2n}((-1)^j-1)\frac{(-1;q)_jq^j}{(q;q)_j}\nonumber\\
&=-2\sum_{n=0}^{\infty}\frac{(q;q^2)_n(-q)^n}{(-q;q^2)_{n}(1+q^{2n})}\sum_{j=1}^{n}\frac{(-1;q)_{2j-1}q^{2j-1}}{(q;q)_{2j-1}}\nonumber\\
&=\frac{4q^2}{(1+q)}\sum_{n=0}^{\infty}\frac{(q^3;q^2)_n(-q)^n}{(-q^3;q^2)_{n}(1+q^{2n+2})}\sum_{j=0}^{n}\frac{(-q;q)_{2j}q^{2j}}{(q^2;q)_{2j}}.
\end{align}}
Invoking \eqref{diff3}, we see that the proof of Theorem \ref{congmod41} is complete.

\end{document}